\documentclass[12pt,a4paper]{article}

\usepackage{amsmath,amsthm,amssymb,amscd,enumerate}
\usepackage[english]{babel}
\usepackage[latin1]{inputenc}
\usepackage{tikz}
\usepackage{listings}
\usepackage{graphicx,graphics}
\usepackage{bbm}
\usepackage{tocbibind}
\usepackage{fancyhdr}
\usepackage{pgf}

\begin{document}

\title{A two-strain ecoepidemic competition model}
\author{Roberto Cavoretto, Simona Collino, Bianca Giardino,\\Ezio Venturino\\
Dipartimento di Matematica ``Giuseppe Peano'',\\
Universit\`a di Torino, \\via Carlo Alberto 10, 10123 Torino, Italy,\\
email: ezio.venturino@unito.it}
\date{}

\maketitle

\begin{abstract}
In this paper we consider a competition system in which two diseases spread by contact.
We characterize the system behavior, establishing that only some configurations are possible.
In particular we discover that coexistence
of the two strains is not possible, under the assumptions of the model.
A number of transcritical bifurcations relate the more relevant system's equilibria.
Bistability is shown between a situation in which only the disease-unaffected population
thrives and another one containing only the second population with endemic disease.
An accurate computation of the separating surface of the basins of attraction
of these two mutually exclusive equilibria is obtained via
novel results in approximation theory.
\end{abstract}

{\textbf{Keywords}: epidemics; competition; disease transmission; ecoepidemics
{\textbf{AMS MR classification} 92D30, 92D25, 92D40

\section{Introduction}

In the last eighty years of the past century a wealth of literature has been devoted to mathematical issues in ecology. Models for interacting
populations have been studied to investigate real life situations that range for instance from the management and conservation of wild populations in reserve
parks to the microscopic level of cellular interactions and proliferation in cancer research, \cite{Castillo-Chavez,MPV,Murray}.
On the other hand, the impact of disease transmission
on human populations is severe. During the same time span, mathematical means have also been developed to assist
epidemiologists in their daily fight against epidemics. The fundamental contribution of mathematical epidemiology to the historical decision of the WHO in 1980
to discontinue on a worldwide basis the vaccination against smallpox officially declaring the disease extinct is not to be underestimated.

Recently, in mathematical epidemiology, more complicated situations than the usual SIRS (Susceptible, Infected, Removed, Susceptible) models
have attracted the attention of researchers.
In nature indeed, it is not unlikely that individuals get infected by more than one disease, in general terms, \cite{CC,Li},
as well as for the case of specific diseases, such as the dengue fever, \cite{EC}. In
particular, tuberculosis has received quite a lot of attention, in view of the fact that its nowadays recrudescence, \cite{HD,JLF,Iannelli}.
But also the more widespread flu has been considered, \cite{ALL,CB}. The rather general situation of multistrain epidemic models is investigated in \cite{AKW}.
There are some
instances in which both strains coexist in the host, and in this case one talks about coinfection models. Alternatively, it may happen
that the most recently acquired disease replaces the older one.
In such case we are in presence of the so-called superinfection phenomenon, \cite{Cai,IML,LLM}.

Population associations, whether they arise for mutual benefit or more generally for the survival of at least one at the expense of others,
are common in nature. In fact mathematical biology research received a great boost from the early works of Lokta and Volterra on predator-prey
systems, \cite{Lo,Vo}. On the other hand,
competition models in ecology are among the first ones investigated, and still arise interest among researchers, \cite{AZ,ADV,CC,W}.
They were also subject of in vitro experiments, recall the well-known Gause laboratory investigations to assess the growth rates of two
bacteria populations, both when living indipendently as well as when they were kept in the same environment. In the last
situation he was able to determine also the rates at which the two populations compete for resources, \cite{Pop_Models}. His results
contributed to support with field data the theoretical results on the logistic model and of competition systems.

The two fields of research described above, namely population theory and mathematical epidemiology,
have almost independently progressed, until the nineties. Then, the first models accounting for diseases spreading by
contact among interacting populations, showing that the demographic equilibria of the systems under consideration were sensibly altered by
the epidemics appeared. A new branch of science was developed, named ecoepidemiology, see Chapter 7 of \cite{MPV} for an account of its early days.
In fact, ecoepidemilogy studies dynamical systems describing populations interactions among which a disease spreads
by contact. The underlying demographic subsystem can be of various types.

Recently, also in the framework of ecoepidemic models, the issue of multiple strains has been considered, \cite{EGV1,EGV2,RRV}.
However, in both these papers, the latter has been investigated only for interactions of predator-prey type. It makes sense to extend
the query to other commonly accepted systems.
In this work therefore we aim at the investigation in the case of competition models.

The paper is organized as follows. We present the model in the next section. Some preliminary results are provided in Section 3.
The disease-free ones are analysed in Section 4 while coexistence with no disease is studied in Section 5.
The following section contains the equilibria where the disease is endemic. Section 7 contains the particular cases of the model,
useful for the final discussion. Bistable situations are investigated in Section 8, before the Conclusions.


\section{The Model}
We introduce here a competition model between two populations in which two recoverable disease strains are present, affecting only one population.
Also, there is no possibility for the individuals to become immune to the diseases.
We further assume that the epidemics are transmitted only horizontally, the diseases are propagated by contact or by demographic
interaction among individuals of the affected population.
We assume that the two diseases affecting one population cannot be transmitted to the other one. Furthermore, we also assume
that the two diseases do not interfere with each other, i.e. there is no superinfection nor coinfection.
This means that at any given time one individual can carry at most one of the
two diseases and in such case cannot be infected by the other one, nor can the second disease replace the former.
We also make the assumption that infected individuals do not reproduce and
that they are too weak also to compete with the other population.
Let $P$ be the first healthy population, while $S$, $V$ and $W$ denote respectively the susceptible individuals of the second population,
the diseased individuals of the first type and those infected by the second strain.
\begin{eqnarray}\label{model}
\dfrac{dP}{dt} = s\left( 1 - \dfrac{P}{L}\right) P - aPS\\ \nonumber
\dfrac{dS}{dt} = r\left( 1 - \dfrac{S}{K}\right) S- bPS - \lambda VS - \beta WS + \psi V + \varphi W \\ \nonumber
\dfrac{dV}{dt} = \lambda VS - \psi V - \mu V -ePV \\ \nonumber
\dfrac{dW}{dt} = \beta WS -\varphi W-\nu W  -fPW
\end{eqnarray} 
In the first equation we describe the dynamics of the disease-free population. It reproduces logistically and is subject to the negative
influence of the competing population for resources, but only the influence of the healthy individuals of the second population is felt.
The second equation contains the healthy second population: again, it reproduces logistically, but in view of the weakeness of the diseased
individuals, their presence is not accounted for in the intraspecific competition term. Instead, the competition with the first population
represents a hindrance for the growth of the second one. Other population losses are due to the infection mechanism, that drives some of their
individuals into the two diseased classes, at different rates since the virulence of the two strains is different. Instead, recovered
individuals migrate back from the infected classes to the healthy one, again at different rates, as the recovery periods for each strain differ.
The last two equations are kind of symmetric, and account for the behavior of infected individuals. They are recruited via successful
contacts between a healthy individual and an infected one, and leave it either by recovery or by natural plus disease-related mortality.
Finally, the pressure of the other population is felt also by individuals in these classes.

In the model (\ref{model}), the parameter $s$ denotes the reproduction rate of the first population,
$L$ its carrying capacity, $a$ the damage inflicted by the susceptible of the second population on the first one;
$r$ the reproduction rate of the second population, $K$ its carrying capacity,
$b$ the damage inflicted by the first population on the susceptibles of the second one, $\lambda$ the first disease contact rate,
$\beta$ the second disease incidence, $\psi$ the first disease recovery rate, $\varphi$ the second disease recovery rate,
$\mu$ the natural plus first strain mortality rate, $e$ the damage inflicted by the first population on the
infected of the first disease, $\nu$ the natural plus second disease mortality rate, $f$ the damage inflicted
by the first population on the infected of the second disease. We assume that all the above parameters are nonnegative.

To study the stability of the equilibria we need the Jacobian matrix of (\ref{model}),
\begin{equation}\label{jac}
J =
\left(
\begin{array}{cccc}
J_{11} & -aP & 0 & 0\\ 
-bS & J_{22} & -\lambda S+\psi & -\beta S +\varphi\\ 
-eV & \lambda V & J_{33} & 0 \\
-fW & \beta W & 0 & J_{44}
\end{array}
\right)
\end{equation}
with
\begin{eqnarray*}
J_{11}=-\frac{sP}{L}+s\left( 1-\frac{P}{L}\right) -aS, \quad J_{22}=-\frac{rS}{K}+\left( 1-\frac{S}{K}\right) -\lambda V -\beta W -bP,\\
J_{33}=\lambda S-\mu - \psi -eP, \quad J_{44}=\beta S-\nu - \varphi - fP.
\end{eqnarray*}

\section{Preliminary results}

In this section we discuss the behaviour of the trajectories when time tends to infinity.\\
 
\subsection{Boundedness}

{\textbf{Theorem.}} The solution trajectories of the model (\ref{model}) are bounded.

{\textbf{Proof.}}
From the model it is easy to show that:
$$
\limsup\limits_{t\to +\infty} P=L, \quad \limsup\limits_{t\to +\infty} S=K.
$$
Thus we assume $P \leqslant L$, $S \leqslant K$. Setting the total environmental population $\Phi (t) = P+S+V+W$, we obtain
\begin{displaymath}
\dfrac{d\Phi}{dt} = s\left( 1-\dfrac{P}{L}\right) P+r\left( 1-\dfrac{S}{K}\right) S-\mu V-\nu W-(a+b)PS-ePV-fPW.
\end{displaymath}
Let be $H:=(a+b)PS+ePV+fPW$ ( $\geqslant 0$ ), we can make the following consideration for all positive $\varepsilon$
\begin{align*}
\dfrac{d\Phi}{dt} +\varepsilon\Phi \\
= s\left( 1-\dfrac{P}{L}\right) P+\varepsilon P+r\left( 1-\dfrac{S}{K}\right) S+\varepsilon S-(\mu-\varepsilon) V-(\nu-\varepsilon) W-H\\
\leqslant s\left( 1-\dfrac{P}{L}\right) L+\varepsilon L+r\left( 1-\dfrac{S}{K}\right) K+\varepsilon K-(\mu-\varepsilon) V-(\nu-\varepsilon) W-H \\
\leqslant sL+\varepsilon L+rK+\varepsilon K-(\mu-\varepsilon) V-(\nu-\varepsilon) W \\
= C-(\mu-\varepsilon) V-(\nu-\varepsilon) W.
\end{align*}
where $C:=sL+\varepsilon L+rK+\varepsilon K$ is a positive constant. Let be $\varepsilon_0 := \min\{\mu,\nu\}$,
for all $\varepsilon$ so that $0<\varepsilon<\varepsilon_0$ we obtain
\begin{displaymath}
\dfrac{d\Phi}{dt} +\varepsilon\Phi \leqslant C
\end{displaymath}
from which
\begin{displaymath}
\Phi (t)\leq \dfrac{C}{\varepsilon}+k e^{-\varepsilon t}\leq M
\end{displaymath}
for some suitable constant $M$, so we have shown the system trajectories are bounded and cannot go to infinity.

\subsection{The possible equilibria}

Of the 16 possible equilibria of (\ref{model}), only 8 are viable. They are summarized in the following Table.

\ \\

\begin{center}
\begin{tabular}{|c|c|c|c|c|}
\hline 
\textbf{P} & \textbf{S} & \textbf{V} & \textbf{W} & \textbf{EQUILIBRIUM} \\ 
\hline 
0 & 0 & 0 & 0 & $E_{0}$ \\ 
\hline 
+ & 0 & 0 & 0 & $E_{1}$ \\ 
\hline 
0 & + & 0 & 0 & $E_{2}$ \\ 
\hline 
+ & + & 0 & 0 & $E_{3}$ \\ 
\hline 
0 & + & + & 0 & $E_{4}$ \\ 
\hline 
0 & + & 0 & + & $E_{5} $ \\ 
\hline 
+ & + & + & 0 & $E_{6} $ \\ 
\hline 
+ & + & 0 & + & $E_{7} $ \\ 
\hline 
\end{tabular}
\end{center}

\section{The disease-free equilibria}

\subsection{Ecosystem preservation}

It is immediate to observe that the origin is a trivial solution of the system (\ref{model}).
The eigenvalues of the Jacobian matrix evaluated at the origin are
$\lambda_{1}=s>0$, $\lambda_{2}=r >0$, $\lambda_{3}= -(\psi+ \mu) <0$ and $\lambda_{4}= -(\varphi + \nu) <0$.
Therefore equilibrium $E_{0}$ is a saddle and so it is unstable. Under our assumptions the system will never be wiped out.


\subsection{Only the population unaffected by the disease survives}
The population unaffected by the disease settles at the environment's carrying capacity level, $P_1=L$, while the other one is completely wiped out.
The eigenvalues of the Jacobian matrix evaluated at the equilibrium are
$\lambda_{1}=-s<0$, $\lambda_{2}=-\mu - \psi -eL <0$, $\lambda_{3}=-\nu -\varphi - fL <0$, $\lambda_{4}=r-bL$.
The stability of the equilibrium depends only on the sign of the eigenvalue $\lambda_{4}$.
The first population survives and settles to its carrying capacity if and only if
\begin{equation}\label{E1_stab}
L>\dfrac{r}{b}.
\end{equation}
If this condition is not satisfied, $E_{1}$ is a saddle.

\subsection{Only healthy population thrives}
Here, it is the disease-affected population that settles at carrying capacity, $S_2=K$, while the disease is eradicated, and also the competitor gets
extinguished.
The eigenvalues of the Jacobian are:
$\lambda_{1}=-r<0$, $\lambda_{2}=s-aK$, $\lambda_{3}=\lambda K -\psi -\mu$, $\lambda_{4}=\beta K -\varphi -\nu$.
Stability depends on the sign of the eigenvalues $\lambda_{2}$, $\lambda_{3}$, $\lambda_{4}$, which
reduce to the following conditions
\begin{equation}\label{E2_stab}
K>\dfrac{s}{a}, \quad K <\dfrac{\psi +\mu}{\lambda}, \quad K < \dfrac{\varphi +\nu}{\beta}.
\end{equation}
If they are not satisfied, $E_{2}$ is a saddle.

\section{Coexistence in a disease-free environment}

The populations levels at this equilibrium are
$$
P_3 = \dfrac{Lr(-s+aK)}{-rs+bLKa}, \quad S_3=\dfrac{Ks(-r+bL)}{-rs+bLKa}.
$$
The equilibrium is feasible if
one of the two alternative conditions holds
\begin{eqnarray}\label{E3_feas1}
aK>s, \quad bL>r ;\\ \label{E3_feas2}
aK<s, \quad bL<r .
\end{eqnarray}
Two eigenvalues of (\ref{jac}) at the equilibrium can be explicitly obtained, $\beta S_3-\varphi-\nu -fP_3$,
$\lambda S_3-\psi-\mu -eP_3$, while the other two are roots of a quadratic.
Explicitly, we find
\begin{eqnarray*}
\lambda_{1}= \dfrac{rs(-bL-aK+r+s)+\sqrt{\Delta}}{2(bLKa-rs)} \\
\lambda_{2}= \dfrac{rs(-bL-aK+r+s)-\sqrt{\Delta}}{2(bLKa-rs)}\\
\lambda_{3}=\dfrac{fLr(s-aK)-\beta Ks(r-bL)}{bLKa-rs}-(\nu +\varphi) \\
\lambda_{4}=\dfrac{eLr(s-aK)-\lambda Ks(r-bL)}{bLKa-rs}-(\mu +\psi)
\end{eqnarray*}
where
\begin{equation}\label{delta}
\Delta :=
rs(rs(bL+aK)^2-2rs(r-s)(bL-aK)+rs(r-s)^2\\
-4bLKa(sLb+raK-bLaK)).
\end{equation}
To attain stability we need all of the eigenvalues to have negative real parts.

All eigenvalues are real if $\Delta \geq 0$, i.e. if
\begin{eqnarray}\label{delta_pos}
rs(bL+aK)^2-2rs(r-s)(bL-aK)+rs(r-s)^2\\ \nonumber
\geq 4bLKa(sLb+raK-bLaK).
\end{eqnarray}
In this case the equilibrium is a node.
There are instead two complex conjugated eigenvalues if $\Delta < 0$, namely if the inequality in (\ref{delta_pos}) is reversed. In such case we have a focus.

The Routh-Hurwitz conditions for stability give in this case of a quadratic equation the following inequalities:
$aS_3+bP_3>0$ which holds always, and $S_3P_3 (KL)^{-1}(rs-abKL)>0$, i.e. $rs>abKL$.
Thus $P_3$ is stable only if the set
of feasibility conditions (\ref{E3_feas2}) are satisfied, and furthermore if
\begin{eqnarray}\label{stab_E3}
s > \max \left\{ \dfrac{LKa(\nu b+\varphi b+fr)}{fLr-\beta Kr+\beta LKb+\nu r+\varphi r},
\dfrac{LKa(\mu b+\psi b+er)}{eLr-\lambda Kr+\lambda LKb+\mu r+\psi r} \right\}.
\end{eqnarray}

\section{Equilibria with endemic disease}
\subsection{First strain endemic equilibrium}

At $E_{4}$ we find the population levels
\begin{align*}
S_4 = \dfrac{\psi +\mu}{\lambda}, \quad V_4=\dfrac{r(\mu+\psi) (K\lambda-\mu -\psi)} {K\lambda^{2}\mu}.
\end{align*}
Recall that $S$ and $V$ represent the number of individuals of the second population respectively susceptible and infected by the first disease.
The feasibility conditions for $E_{4}$ are
\begin{equation}\label{E4_feas}
K\ge \dfrac{\psi +\mu}{\lambda} =: A.
\end{equation}
If we set
\begin{align*}
\Delta:=r^2 \mu^4+2r^2 \mu^2 K\lambda\psi-2r^2 \mu^2 \psi^2+r^2 K^2 \lambda^2 \psi^2 -2r^2 K\lambda\psi^3 + r^2 \psi^4\\
-4K^2 \mu^3 r\lambda^2 -4K^2 \mu^2 r\lambda^2 \psi+4K\mu^2 r\lambda\psi^2+4K\mu^4 r\lambda+8K\mu^3 r\lambda\psi
\end{align*}
we obtain the following eigenvalues:
\begin{eqnarray*}
\lambda_{1}= \dfrac{-rC+\sqrt{\Delta}}{2K\lambda\mu}, \quad
\lambda_{2}= \dfrac{-rC-\sqrt{\Delta}}{2K\lambda\mu},\quad
\lambda_{3}=\beta A-\nu -\varphi, \quad
\lambda_{4}=s-aA,
\end{eqnarray*}
where we have introduced new parameters as follows
$$
B:=\dfrac{\varphi+\nu}{\beta},\quad
C:=\mu^2+K\lambda\psi-\psi^2.
$$

All eigenvalues are real if $\Delta \geq 0$, i.e. for
\begin{align*}
r\mu^2(\mu^2 +2K\lambda\psi)+r\psi^2(K^2\lambda^2 -2\mu^2)+4K\mu^2\lambda(\mu+\psi)(\mu+\psi-K\lambda) \geq r\psi^3(2K\lambda-\psi)
\end{align*}

The stability conditions are
$$
\sqrt{\Delta} < rC, \quad
\sqrt{\Delta} > -rC, \quad
\dfrac{s}{a}<A<B.
$$
The first two inequalities coexist if and only if
$-rC<\sqrt{\Delta}<rC$ i.e. if and only if $C>0$, which is true in view of (\ref{E4_feas}). In fact it implies
$$
K>\dfrac{\psi^2-\mu^2}{\lambda\psi}=\dfrac{\psi-\mu}{\psi}A
$$
and the last term is obviously smaller than $A$, so that if $E_4$ is feasible, the above condition holds.
In conclusion $E_{4}$ is a stable node if and only if
\begin{equation}\label{E4_stab_n}
\dfrac{s}{a}<A<B, \quad
\sqrt{\Delta}<rC.
\end{equation}

There are instead two complex conjugate eigenvalues if $\Delta <0$ and 
in this situation to attain stability we need $r(\mu^2+K\lambda\psi^2-\psi^2)>0$ which becomes
$$
K>\dfrac{\psi-\mu}{\psi}A .
$$
This request is always true in view of (\ref{E4_feas}).

Stability of the focus is thus guaranteed by
\begin{equation}\label{E4_stab_f}
\frac{s}{a}<A<B.
\end{equation}

\subsection{Second strain endemic equilibrium}

The equilibrium $E_5$ has the population levels
\begin{align*}
S_5 = \dfrac{\varphi +\nu}{\beta}, \quad W_5=\dfrac {r(\varphi +\nu)(K \beta -\varphi-\nu)}{K\beta^{2}\nu}
\end{align*}
The feasibility conditions for $E_{5}$ become
\begin{equation}\label{E5_feas}
K>\dfrac{\varphi +\nu}{\beta} =: \widetilde B
\end{equation}
If we set 
\begin{align*}
\widetilde {\Delta}:=r^2 \nu^4+2r^2 \nu^2 K\beta\varphi-2r^2 \nu^2 \varphi^2+r^2 K^2 \beta^2 \varphi^2 -2r^2 K\beta\varphi^3 + r^2 \varphi^4\\
-4K^2 \nu^3 r\beta^2 -4K^2 \nu^2 r\beta^2 \varphi+4K\nu^2 r\beta\varphi^2+4K\nu^4 r\beta+8K\nu^3 r\beta\varphi
\end{align*}
as well as
$$
\widetilde A:=\dfrac{\psi+\mu}{\lambda}, \quad
\widetilde B:=\dfrac{\varphi+\nu}{\beta}, \quad
\widetilde D:=\nu^2+K\beta\varphi-\varphi^2,
$$
we obtain the following eigenvalues
$$
\lambda_{1}= \dfrac{-r\widetilde D+\sqrt{{\widetilde {\Delta}}}}{2K\beta\nu}, \quad
\lambda_{2}= \dfrac{-r\widetilde D-\sqrt{{\widetilde {\Delta}}}}{2K\beta\nu}, \quad
\lambda_{3}=\lambda \widetilde B-(\mu +\psi), \quad
\lambda_{4}=s-a\widetilde B.
$$
Assuming that ${\widetilde {\Delta}} \geq 0$, 
we find the following conditions of stability:
$$
-r\widetilde D < \sqrt{\widetilde {\Delta}} < r\widetilde D, \quad
\dfrac{s}{a}< \widetilde B<\widetilde A.
$$
The first inequalities hold if and only if ${\widetilde {D}} > 0$.
But from the feasibility condition $K>\widetilde B$ this condition in fact follows. Indeed it is equivalent to
$\nu^2+K\beta\varphi-\varphi^2 >0$ which can be rewritten as
$$
K>\dfrac{\varphi^2-\nu^2}{\beta\varphi}=\dfrac{\varphi-\nu}{\varphi}\widetilde B.
$$
The last quantity is smaller than $\widetilde B$, so that if (\ref{E5_feas}) holds, the claim follows.

In conclusion $E_{5}$ is a stable node if and only if
\begin{equation}\label{E5_stab_n}
\dfrac{s}{a}<\widetilde B<\widetilde A, \quad
\sqrt{\widetilde {\Delta}}<r\widetilde D.
\end{equation}

The eigenvalues are instead complex conjugate if $\widetilde {\Delta} <0$ or rather if
\begin{align*}
r\nu^2(\nu^2 +2K\beta\varphi)+r\varphi^2(K^2\beta^2 -2\nu^2)+r\varphi^3(\varphi-2K\beta)+4K\nu^2\beta(\nu+\varphi)(\nu+\varphi-K\beta) < 0 \\
\end{align*}
and for stability we need
$r(\nu^2+K\beta\varphi^2-\varphi^2)>0$ which can be rewritten as 
$$
K>\dfrac{\varphi-\nu}{\varphi}\widetilde B,
$$
which always holds as above in view of (\ref{E5_feas}).
$E_5$ is a stable focus if
\begin{equation}\label{E5_stab_f}
\dfrac{s}{a}<\widetilde B<\widetilde A
\end{equation}

\subsection{Coexistence of first strain
and healthy population.}

At $E_{6}$ we find
$$
P_6 = -\dfrac{L(-\lambda s+a\mu +a\psi)}{\lambda s+eLa}, \quad
S_6 = \dfrac{s(\mu+\psi+eL)}{\lambda s+eLa}, \quad
V_6=\dfrac{s(\mu+\psi+eL)E}{K(\lambda s+eLa)F},
$$
where we have defined
\begin{eqnarray*}
E:=rLKea-eLrs-LKb\lambda s+bLKa\psi+bLKa\mu-rs\psi+rK\lambda s-rs\mu \\
F:=s\lambda eL+s\lambda\mu -\psi eLa
\end{eqnarray*}
For the population $P$ to be nonnegative we require $-\lambda s+a\mu +a\psi<0$ i.e.
\begin{equation}\label{P6_pos}
a<\dfrac{\lambda s}{\mu+\psi}.
\end{equation}
For the nonnegativity of $V$ we need either one of the following pairs of conditions
\begin{align*}
E>0, \quad F>0; \qquad E<0, \quad F<0.
\end{align*}
Introducing new parameters as follows,
$$
M:=\dfrac{s(eLr+bLK\lambda+r\psi-rK\lambda+r\mu)}{LK(er+b\psi+b\mu)}, \quad
N:=\dfrac{s\lambda(eL+\mu)}{eL\psi}, \quad
G:=\dfrac{\lambda s}{\mu+\psi}
$$
to rewrite the nonnegativity of $V$, the feasibility conditions can be written as the two alternatives
\begin{equation}\label{E6_feas_0}
M<a<N, \quad a<G; \qquad N<a<M,\quad a<G.
\end{equation}
But the second case is impossible, since $G<N$. In fact, this inequality explicitly is
$$
\dfrac{\lambda s}{\mu+\psi} < \dfrac{s\lambda(eL+\mu)}{eL\psi}\\
$$
which becomes
$$
\dfrac{-eL\mu-\mu^2-\mu\psi}{eL\psi(\mu+\psi)}<0
$$
and the latter holds unconditionally.

From the first case of (\ref{E6_feas_0}), using the previous result, instead we obtain the feasibility condition for $E_6$
\begin{equation}\label{E6_feas}
M<a<G.
\end{equation}
Note that $M<G$ gives
\begin{displaymath}
K\geq \dfrac{\mu+\psi}{\lambda}.
\end{displaymath}

The search of the condition of stability is too complicated both studying the sign of the eigenvalues and applying the criterion of Routh-Hurwitz.

\subsection{Second strain
with healthy population coexistence.}

The Equilibrium $E_{7}$ has the population values
\begin{eqnarray*}
P = -\dfrac{L(-\beta s+a\nu +a\varphi)}{\beta s+fLa}, \quad
S = \dfrac{s(\nu+\varphi+fL)}{\beta s+fLa}, \quad
W=\dfrac{s(\nu+\varphi+fL)\widehat E}{K(\beta s+fLa)\widehat F}
\end{eqnarray*}
where
\begin{eqnarray*}
\widehat E:=rLKfa-fLrs-bLK\beta s+bLKa\varphi+bLKa\nu-rs\varphi+rK\beta s-rs\nu, \\
\widehat F:=s\beta fL+s\beta\nu -\varphi fLa.
\end{eqnarray*}
$P$ is nonnegative if and only if
$$
a<\dfrac{\beta s}{\nu+\varphi},
$$
while $W$ is nonnegative if and only if either one of the two sets of inequalities holds
\begin{align*}
\widehat E>0, \quad
\widehat F>0 ; \qquad
\widehat  E<0, \quad
\widehat F<0
\end{align*}
To explicitly study these conditions we introduce the new parameters:
$$
\widehat M:=\dfrac{s(fLr+bLK\beta+r\varphi-rK\beta+r\nu)}{LK(fr+b\varphi+b\nu)}, \quad
\widehat N:=\dfrac{s\beta(fL+\nu)}{fL\varphi}, \quad
\widehat G:=\dfrac{\beta s}{\nu+\varphi},
$$
so that they become the two alternative sets
\begin{align*}
\widehat M<a<\widehat N ,\quad
a<\widehat G ; \qquad
\widehat N<a<\widehat M, \quad
a<\widehat G
\end{align*}
Again, the second case is impossible, because $G<N$. In fact
we have the inequality
$$
\dfrac{\beta s}{\nu+\varphi} < \dfrac{s\beta(fL+\nu)}{fL\varphi}
$$
from which
$$
\dfrac{-fL\nu-\nu^2-\nu\varphi}{fL\varphi(\nu+\varphi)}<0
$$
which is always true.

The first case reduces to
\begin{equation}\label{E7_feas}
\widehat M<a<\widehat G
\end{equation}
giving the feasibility condition for $E_{7}$,
which explicitly can be written as
\begin{displaymath}
K\geq \dfrac{\nu+\varphi}{\beta}
\end{displaymath}

\section{Analysis of the particular cases}
\subsection{The diseases-free system}

In this section we study the behaviour of the demographic competition model underlying the ecoepidemic system,
to compare the results with those previously obtained, so as to highlight the effect of the diseaeses.
Omitting the epidemics, the model is reduced to
$$
\dfrac{dP}{dt} = s\left( 1 - \dfrac{P}{L}\right) P - aPS, \quad
\dfrac{dS}{dt} = r\left( 1 - \dfrac{S}{K}\right) S- bPS,
$$
whose equilibria are
$Q_{0}=(0,0)$, $Q_{1}=(L,0)$, $Q_{2}=(0,K)$ and 
$$
Q_{3}=\left( \dfrac{Lr(aK-s)}{bLKa-rs},\dfrac{Ks(bL-r)}{bLKa-rs} \right)
$$
The Jacobian $\widetilde J$ is 
$$
\left(
\begin{array}{cc}
-\dfrac{sP}{L}+s\left( 1-\dfrac{P}{L}\right) -aS & -aP \\ 
-bS & -\dfrac{rS}{K}+r \left( 1-\dfrac{S}{K}\right) -bP  
\end{array}
\right)
$$

It is immediately seen that the origin is unstable, given the eigenvalues $s$ and $r$.
At $Q_1$, whose nonzero components are those of $E_1$, we find instead 
$-s$ and $r-bL$, so that the stability of the equilibrium hinges on the very same condition (\ref{E1_stab}) for $E_1$.

For $Q_2$, which coincides with the nontrivial part of $E_2$, we have 
$-r$ and $s-aK$, so that stability is ensured by
\begin{equation}\label{Q2_stab}
K>\dfrac{s}{a}.
\end{equation}

Feasibility of $Q_{3}$ is the very same (\ref{E3_feas1}) or  (\ref{E3_feas2})
as for $E_3$ of which the former represents the projection onto the $S-P$ phase subspace.

The Routh-Hurwitz conditions for stability give
$$
-{\textrm{tr}} (\widetilde J)=\frac sL P_3 + \frac rK S_3>0
$$
which is clearly satisfied, and
$$
\det (\widetilde J)= \frac {S_3P_3}{LK} (rs-abKL)>0.
$$
This condition is clearly incompatible with (\ref{Q2_stab}) and (\ref{E1_stab}). Therefore whenever the mutually
exclusive equilibria $Q_1$ and $Q_2$ are both stable, $Q_3$ is unstable and vice versa.

\subsection{The one disease epidemic model}

Omitting the $W$ strain, it is easily seen that the equilibria of the subsystem of (\ref{model}) in which $P$ and $W$ are absent
are the origin, which is unstable, the disease-free point $(K,0)$, which is stable when the second
(\ref{E2_stab}) is satisfied, and the endemic equilibrium $(S_*,V_*)\equiv (S_4,V_4)$. The latter is feasible
when (\ref{E4_feas}) holds. In view of these conditions, there is a transcritical bifurcation when both the second of
(\ref{E2_stab}) and (\ref{E4_feas}) become equalities from which the equilibrium containing the infected subpopulation
emanates from the disease-free equilibrium
and the disease establishes itself endemically in the system.

Evidently, similar conclusions hold for the $W$ strain in absence of the $V$-affected individuals.

\section{Bistability}

Bistability is achieved for the following set of parameters
$s = 0.4$, $L = 1.5$, $a = 0.3$, $b = 0.7$, $e = 0.2$, $f = 0.2$, $r = 0.7$, $K = 2$, $\lambda = 0.7$, $\psi = 0.2$, $\mu = 0.5$,
$\varphi=0.7$, $\nu=0.9$, $\beta=0.2$.
Taking the initial condition as
$(0.0,1.8,0.1,0.1)$
we obtain the healthy population-free equilibrium with endemic disease in the second population
$E_4=(0,S,V,0)=(0, 1, 0.7, 0)$, see Figure \ref{fig:healthy-free},
while taking the point $(1.7,0.8,0.1,0.1)$ we find the diseased population-free equilibrium
$E_1=(P,0,0,0)=(1.5, 0, 0, 0)$, see Figure \ref{fig:disease-free}.
Instead, allowing for a nonzero initial value of the $P$ population, namely $(0.1,1.8,0.1,0.1)$,
we obtain the $E_6=(P,S,V,0)=(0.2828, 1.0760, 0.2441, 0)$ equilibrium, see Figure \ref{fig:coexistence}.

\begin{figure}[htb]
\centering
\includegraphics[scale=0.8]{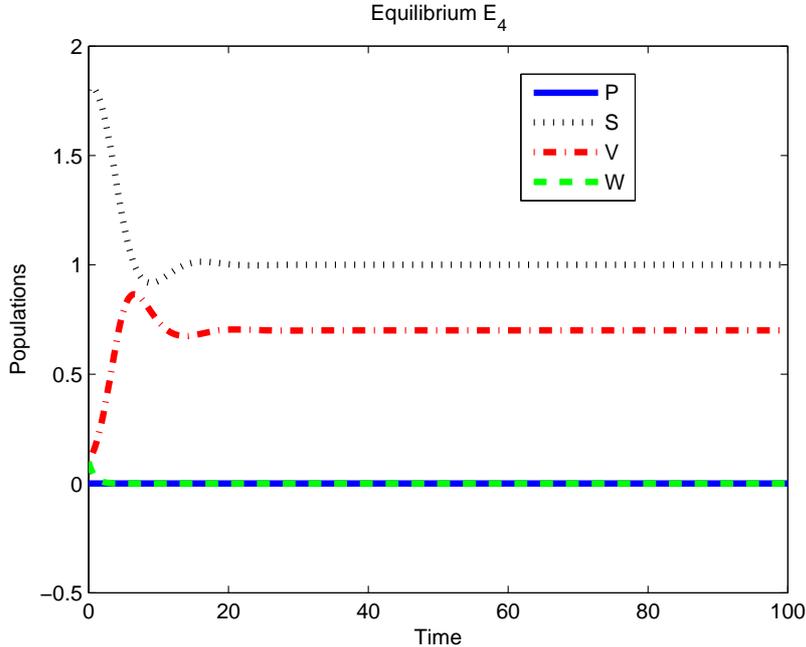}
\caption{The healthy population and second disease-free equilibrium $E_4=(0,S,V,0)=(0, 1, 0.7, 0)$, obtained with the parameter values
$s = 0.4$, $L = 1.5$, $a = 0.3$, $b = 0.7$, $e = 0.2$, $f = 0.2$, $r = 0.7$, $K = 2$, $\lambda = 0.7$, $\psi = 0.2$, $\mu = 0.5$,
$\varphi=0.7$, $\nu=0.9$, $\beta=0.2$
and initial conditions $(0.0,1.8,0.1,0.1)$. Here the $P$ and $W$ lines overlap and only one is shown.}
\label{fig:healthy-free}
\end{figure}

\begin{figure}[htb]
\centering
\includegraphics[scale=0.8]{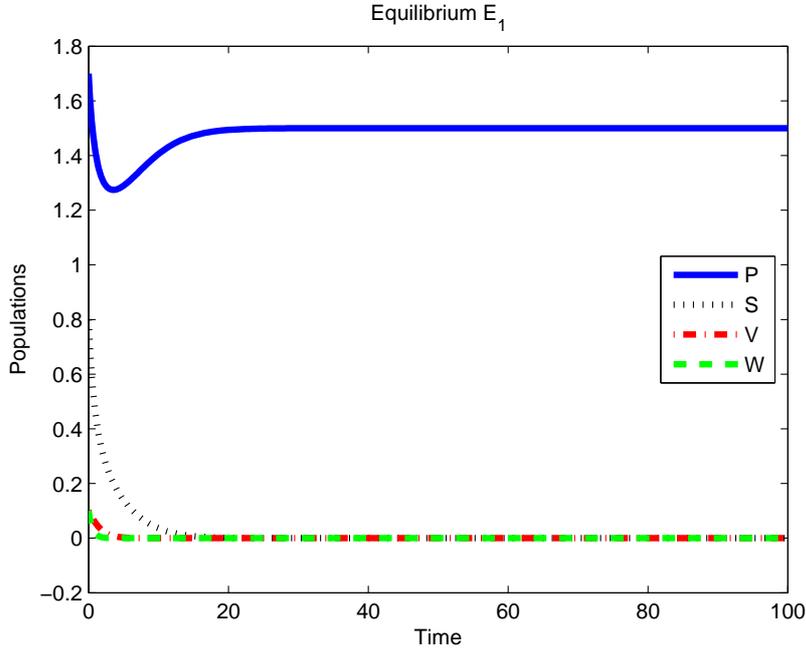}
\caption{The diseased-free equilibrium $E_1=(P,0,0,0)=(1.5, 0, 0, 0)$, obtained with the same parameter values of Figure \ref{fig:healthy-free}
and initial conditions $(1.7,0.8,0.1,0.1)$.}
\label{fig:disease-free}
\end{figure}

\begin{figure}[htb]
\centering
\includegraphics[scale=0.8]{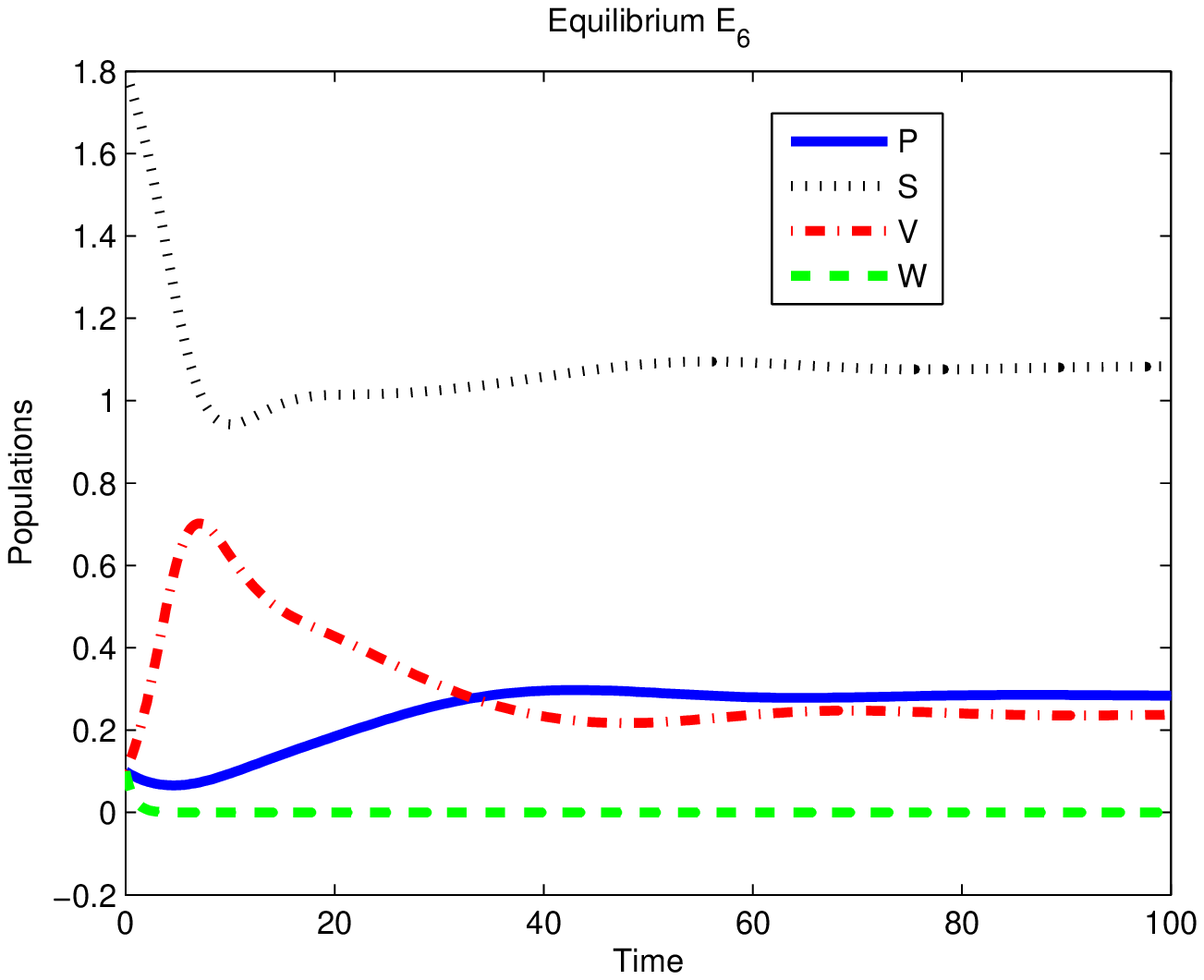}
\caption{The diseased population-free equilibrium $E_6=(P,S,V,0)=(0.2828, 1.0760, 0.2441, 0)$,
obtained with the same parameter values of Figure \ref{fig:healthy-free}
and initial conditions $(0.1,1.8,0.1,0.1)$.}
\label{fig:coexistence}
\end{figure}

In the three-dimensional $P-S-V$ phase space projection of the four dimensional phase space $P-S-V-W$,
it appears valuable here to explicitly assess the surface separating the basins of attractions of these two equilibria. This is achieved
via an algoritm described in \cite{ICNAAM,CMMSE}, see Figure \ref{fig:separatrix}.
The parameters used are $s = 0.3$, $L = 1.5$, $r = 0.7$, $K = 3$, $\lambda = 0.6$,
$\psi = 0.8$, $\mu = 0.3$, $a = 0.2$, $b = 0.5$, $e = 0.2$.

\begin{figure}[htb]
\centering
\includegraphics[scale=0.4]{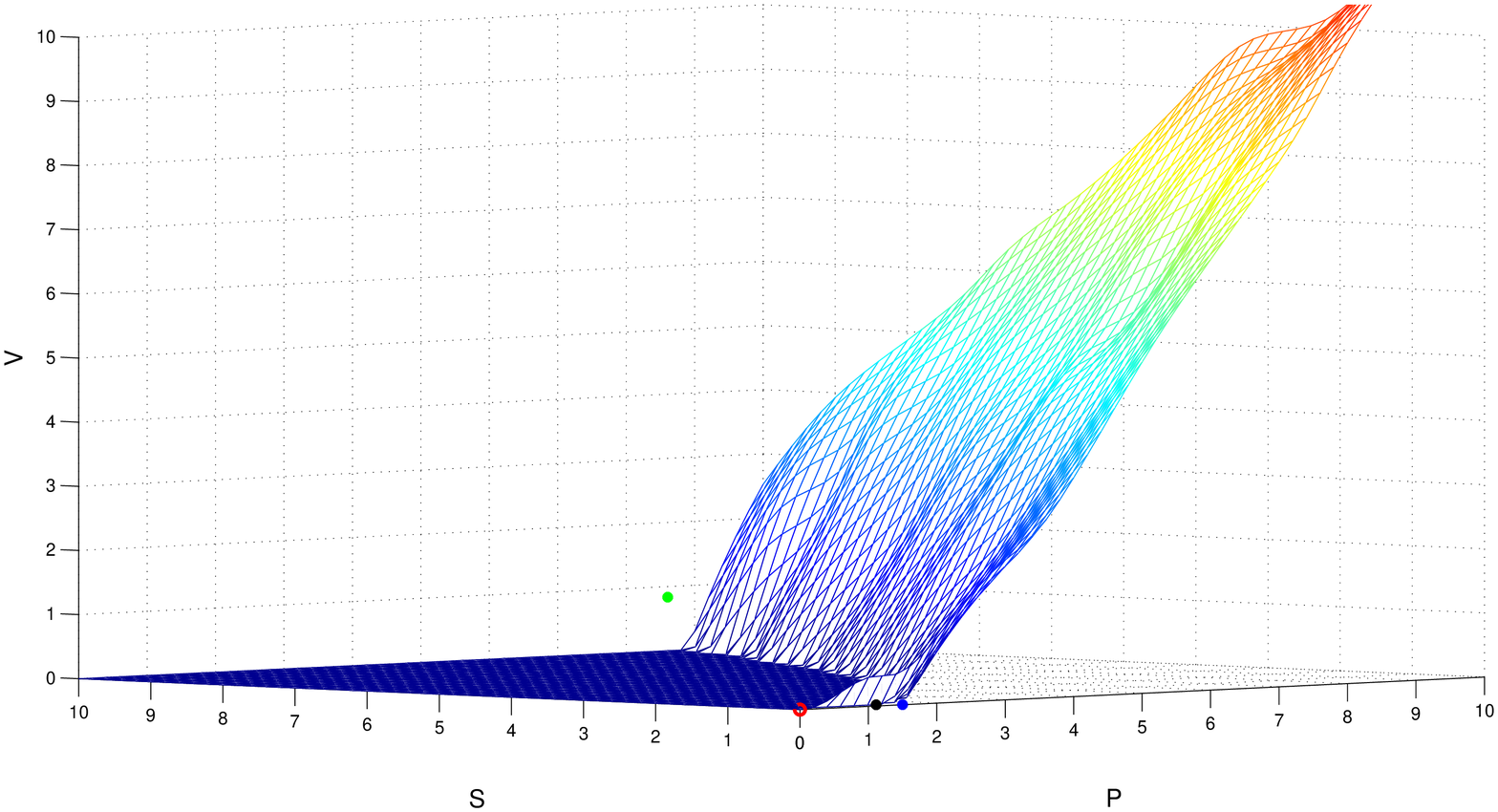}
\caption{The separatrix in the $W=0$ phase subspace
obtained with the parameters $s = 0.3$, $L = 1.5$, $r = 0.7$, $K = 3$, $\lambda = 0.6$,
$\psi = 0.8$, $\mu = 0.3$, $a = 0.2$, $b = 0.5$, $e = 0.2$. The origin is the red dot on the front of the figure,
the green point on the left of the surface is equilibrium $E_4$, the blue point on the axis on the right of
the surface is equilibrium $E_1$, the black point on the surface is the saddle $E_3$.
}
\label{fig:separatrix}
\end{figure}

\section{Conclusions}
In this work we presented a model of competition between two populations characterized by two disease
strains affecting only one of them. In particular this investigation differs in the underlying demographics
model from the systems considered in
\cite{EGV1,EGV2} in that the latter papers consider predator-prey models, with diseases in the prey,
and from \cite{RRV}, where the two epidemics affect the predators.

The model analysis indicates that its trajectories are ultimately bounded and further it
states the presence of seven possible equilibrium points,
since the system's collapse is shown to be impossible.
The rather surprising result is that
no equilibrium allows coexistence of all the four subpopulations. This parallels the results of
the other former predator-prey ecoepidemic model investigations, both in the case of the disease affecting the prey, \cite{EGV1,EGV2},
as well as the predators, \cite{RRV}.

Evidently, from the ecological point of view of biodiversity and for epidemiological considerations,
the best equilibrium that the system can achieve is the coexistence of
the two healthy populations, $E_3$. In a pure competition model in general one knows that the principle of competitive
exclusion holds, but coexistence $Q_3$ is nevertheless also possible. We have found that the same occurs also in the
two-strain ecoepidemic model. In fact, conditions (\ref{E1_stab}) and (\ref{E2_stab}) can both hold for the very same
choice of parameter values, indicating then bistability, i.e. the mutual exclusion of these two possibilities.
Indeed for the parameters
$s = 0.4$, $L = 0.5$, $a = 0.3$, $r = 0.7$, $K = 1$, $b = 0.7$, $\lambda = 0.7$, $\beta=0.2$, $\psi = 0.2$,
$\varphi=0.7$, $\mu = 0.5$, $\nu=0.9$, $e = 0.2$, $f=0.2$, $E_3$ is achieved, see Figure \ref{fig:E3}.

\begin{figure}[htb]
\centering
\includegraphics[scale=0.8]{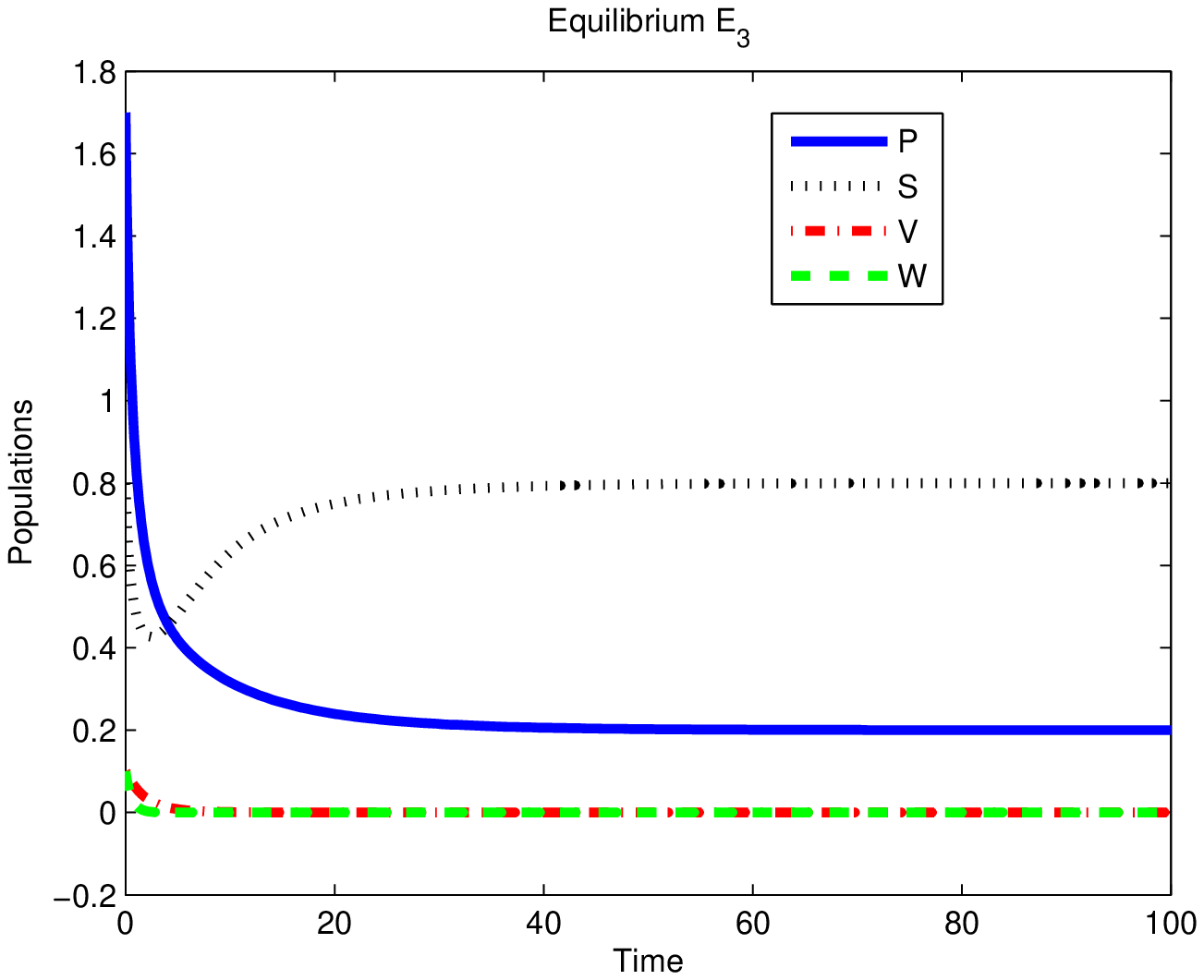}
\caption{Coexistence in a disease-free environment is possible also in the ecoepidemic model,
with the following parameter values
$s = 0.4$, $L = 0.5$, $a = 0.3$, $r = 0.7$, $K = 1$, $b = 0.7$, $\lambda = 0.7$, $\beta=0.2$, $\psi = 0.2$,
$\varphi=0.7$, $\mu = 0.5$, $\nu=0.9$, $e = 0.2$, $f=0.2$.}
\label{fig:E3}
\end{figure}

Stability of the disease-unaffected-population-only equilibrium $E_1$ and its demographic counterpart $Q_1$ coincide.
The stability of the healthy-individuals-only of the infected population equilibrium, $E_{2}$,
depends instead on more conditions than the same equilibrium in the disease-free model, $Q_{2}$.
These conditions involve the epidemics parameters. Therefore, the competitive exclusion principle does not immediately
transfer to the ecoepidemic situation, in that $E_1$ coexists with equilibria
other than $E_2$, as shown by the bistability examples provided.
These other points contain endemically one of the disease strains.

The transcritical bifurcation found for the classical epidemic model has also counterparts in the ecoepidemic system,
because from equilibrium $E_2$ we can see that both $E_4$ and $E_5$ emanate.
Compare indeed the stability conditions for the former,
(\ref{E2_stab}) with the feasibility conditions of both the latter points, (\ref{E4_feas}) and (\ref{E5_feas}).

Again, transcritical bifurcations arise between $E_4$ and $E_6$, when also the first population establishes itself in the ecosystem.
Indeed, recalling the definition of $G$ and $A$, see (\ref{E6_feas_0}) and (\ref{E4_feas}), we see that the second inequality in
(\ref{E6_feas}) holds whenever the first inequality in both
(\ref{E4_stab_f}) and (\ref{E4_stab_n}) is violated, and vice versa. Similar results hold for the second strain,
namely for $E_5$ and $E_7$.

The occurrence of several possible bistability situations with radically differing mutually exclusive equilibria stresses the importance
of the accurate assessment of their basins of attraction. We have provided a step in that direction, with the accurate numerical
determination of the separatrix surface using a novel algoritm explicitly designed for this purpose.

In summary, as it happens for the now standard ecoepidemic models examined usually in the literature, in the case of food chains as well,
the diseases affect heavily the dynamics of the underlying demographic systems. They must therefore be included in the modelling
efforts of theoretical ecologists, in order to arrive at a more accurate description of the natural situations that are being investigated
and thus eventually obtain more reliable results for the policies to be employed for the management of ecosystems.

\end{document}